\documentclass[12pt]{article}
\usepackage{amssymb}
\makeatletter 
\date{}

\newtheorem{theorem}{Theorem}[section]

\newtheorem{proposition}[theorem]{Proposition}

\newtheorem{lemma}[theorem]{Lemma}
\newtheorem{question}[theorem]{Question}
\newtheorem{corollary}[theorem]{Corollary}

\newtheorem{fact}[theorem]{Fact}

\newcommand{\z}{{\Bbb Z}}

\newcommand{\Int}{{\rm Int}}

\newcommand{\lo}{\longrightarrow}

\newcommand{\black}{{\blacksquare}}

\newcommand{\conv}{\overline{\rm conv}}

\begin{document}

\title{COVERING $L^p$ SPACES BY BALLS}
\author{Vladimir P. Fonf
\footnote{Research of the first author was supported in part
by Israel Science Foundation, Grant \# 209/09 and by the Istituto
Nazionale di Alta Matematica of Italy.}, 
Michael Levin
\footnote{Research of
the second author was supported in part by Israel Science Foundation,  Grant \# 836/08}
 and Clemente Zanco
 \footnote{Research of
the third author was supported in part by the Ministero
dell'Universit\`a e della Ricerca Scientifica e Tecnologica of
Italy and by the Center for Advanced Studies in Mathematics at the
Ben-Gurion University of the Negev, Beer-Sheva, Israel.}}

\maketitle

\bigskip
\bigskip
\bigskip

\begin{abstract}
We prove that, given any covering of any separable infinite-dimensional uniformly rotund and uniformly smooth Banach space $X$ by closed balls each of positive radius, some point exists in $X$ which belongs to infinitely many balls. 
\end{abstract}

\bigskip
\bigskip

\bigskip


\bigskip


\vfill\eject

\section{Introduction and statement of the main result}

In the present paper $X$ always denotes a real separable infinite-dimensional Banach space; by {\it ball} in $X$ we mean a closed ball.  Let $\cal A$ be a collection of subsets of $X$.
We say that $\cal A$ is {\it point-finite} if
every point of $X$ belongs to at most finitely many elements of $\cal A$.
A point $x$ of $X$ is said to be a {\it regular point} for $\cal A$ if
there is a neighborhood of $x$ that meets at most finitely many elements of $\cal A$; 
 $x$ is said to be a {\it singular point} for $\cal A$ otherwise. If every point of $X$ is a regular point for $\cal A$,
we say that $\cal A$ is {\it locally-finite} (clearly, that is equivalent to the requirement that no compact set in $X$ meets
infinitely many members of $\cal A$). $\cal A$ is said a {\it covering of X} if each point of $X$ belongs to some member of $\cal A$. 
 
The aim of the present paper is to take a step forward to answering the following 

\begin{question} 
\label{question}
\medskip
{\it Which infinite-dimensional Banach spaces admit point-finite coverings by balls (each of positive radius)?}
\medskip
\end{question}

In order to explain how such a question arises, we recall the following two results (the first one is well known
as ``the Corson's Theorem'').

\begin{theorem} {\rm (\cite{Cor})}
\label{corson}
No (infinite-dimensional) reflexive Banach space admits 
locally finite coverings by bounded closed convex sets.
\end{theorem}

\begin{theorem} {\rm (\cite{MZ})}
\label{index2}
Any real Banach space $X$ can be covered by bounded closed convex sets, each with nonempty interior,
in such a way that no point of $X$ belongs to more than two of them.
\end{theorem}

The family of sets exhibited in the general construction used to prove Theorem \ref{index2} is very far from being a family of balls in the original norm of $X$. Moreover, some classical Banach spaces admit point-finite coverings by balls. For instance, it is easy to check that the covering of $c_0$ that can be obtained by translating the unit ball without overlapping interiors is even locally finite. V. Klee proved in \cite{Kl} that the space $l_1(\Gamma)$ for suitable (uncountable) $\Gamma$ can be covered by translates of its unit ball without overlapping them at all. So Question \ref{question} seems to be very natural, though providing it with a complete answer seems not to be an  easy matter. A first step in that direction have been recently made with the following Theorem, that excludes Hilbert spaces from the class of spaces Question \ref{question} asks for. Even if, looking for spaces outside that class, Hilbert spaces appear as the simplest ones to be considered, up to now no elementary argument for getting such exclusion seems to be available. 

\begin{theorem} {\rm (\cite{FZ}, Theorem 3.2)}
\label{Hilbert}
No covering by balls, each of positive radius, of the infinite-dimensional separable Hilbert space can be point-finite.
\end{theorem}

We refer to \cite{FZ} also for more details and references concerning the subjects involved in the present Introduction.
The goal of this paper is to extend Theorem \ref{Hilbert} to a considerably wider class of spaces. In fact we prove the following

\begin{theorem} {\rm (Main result)} 
\label{goal}
No covering by balls, each of positive radius, of any infinite-dimensional uniformly 
rotund and uniformly smooth separable Banach space can be point-finite.
\end{theorem}

It is very well known that those Banach spaces that are uniformly rotund or uniformly smooth are reflexive
(in fact, super-reflexive). Moreover, if a Banach space $X$ is uniformly rotund or uniformly smooth, then an equivalent norm can be put on $X$ under which $X$ is both uniformly rotund and uniformly smooth. Among those spaces that are both uniformly rotund and uniformly smooth there are $L^p(\mu)$ spaces for any measure $\mu$ and $p \in (1, +\infty)$, so in particular we claim

\begin{corollary} 
\label{L^p}
No covering by balls, each of positive radius, of a separable infinite-dimensional $L^p(\mu)$ space, $1<p<\infty$, $\mu$ any measure, can be point-finite.
\end{corollary}
The proof of Theorem \ref{goal} we provide here is based on a key result of \cite{FZ}; however, after that, it follows a completely different way than what was used in \cite{FZ} to get Theorem \ref{Hilbert}. In fact the argument there runs as follows. Separable polyhedral Banach spaces are first characterized as those whose unit sphere under some equivalent norm admits a {\it point-finite} covering by slices of the unit ball that do not contain the origin. (Recall that a
Banach space is called ``polyhedral'' if the unit ball of any its finite-dimensional subspace is a polytope.)  As a consequence, 
if the unit sphere of some separable Banach space $X$ admits such a covering, then $X$ must be isomorphically polyhedral. It is well known that no (infinite-dimensional) Hilbert space is isomorphically polyhedral. Next, point-finite coverings
of the Hilbert space by balls (if any) are easily reduced to point-finite coverings of the unit sphere by balls that do not contain the origin. Finally, to get a contradiction, these coverings are reduced to point-finite coverings of the unit sphere by slices of the unit ball via the following observation: whenever two spheres in an inner product space do not coincide and have nonempty intersection, such an intersection lies in some hyperplane; this hyperplane splits each of the two balls determined by those spheres in two complementary slices. Unfortunately such a situation characterizes inner product spaces (see \cite{Ami} (15.17)), so the argument cannot be applied outside that class of spaces. 

Our argument here in proving Theorem \ref{goal} has an essential topological component: Corson's Theorem \ref{corson}, which is based on Brouwer's fixed point Theorem, is now our basic tool. We use it in connection
with suitable considerations of geometrical nature.

\medskip

Throughout the paper 
we use standard Geometry of Banach Spaces
notation as in \cite{JL}. In particular, for $x \in X$ and $r>0$, $B(x,r)$ and $S(x,r)$ respectively denote the closed ball and the sphere with center at $x$ and radius $r$; moreover, $B(0,1)$ and $S(0,1)$ are denoted in short respectively by $B_X$ and $S_X$.

\section{Proof of the main result}

In order to prove Theorem \ref{goal},
first of all we notice that, $X$ being separable, we can confine ourselves to prove our theorem for countable coverings. (In fact, let $\{x_n\}$ be any sequence dense in $X$: since each ball has nonempty interior, for some $n_0$ it must happen that $x_{n_0}$ belongs to uncountably many balls.)

We start by borrowing from \cite{FZ} the following Proposition. It describes
a quite general situation in which a sequence of slices of the unit ball 
of any separable Banach space cannot be point-finite.

\bigskip

\begin{proposition} {\rm (\cite{FZ}, Proposition 2.1)}
\label{old}
Let $\{f_i\}_{i=1}^\infty$ be a sequence of norm-one linear functionals on $X$ and $\{\alpha_i\}_{i=1}^\infty$ a sequence in the interval $(0,1)$ converging to $0$. Then the sequence $\{S_i\}_{i=1}^\infty$ of slices of $B_X$ defined by
$$S_i = \{x \in B_X: f_i(x) \geq \alpha_i \}, \ i = 1,2,... $$
is not point-finite.
\end{proposition}

\bigskip

Next we point out a very simple (probably well known) fact.

Roughly speaking, it simply states that any point of a sphere of any uniformly smooth Banach space admits
"almost flat" neighborhoods (relative to the sphere) of "big" diamater provided that the radius of the sphere is "big".
We make this sentence precise in the following way.

\begin{fact}
\label{F1}
Let $X$ be uniformly smooth. For any $\varepsilon > 0$ there exists $b >0$ such that,  
for any $R > b$ and $x \in RS_X$, if $\Gamma_x$ is the hyperplane supporting $RB_X$ at $x$, then
\begin{equation} 
\label{US}
{\rm dist}(y, \Gamma_x) \leq \varepsilon \ \ \ \forall y \in RS_X \cap B(x,2).
\end{equation}

\end{fact}

{\bf Proof.} 
For $t \in S_X$, denote by $f_t$ the (only) norm-one linear functional such that $f_t(t)=1$. Fix $\varepsilon >0$.
By definition of uniform smoothness, for any $\varepsilon >0$ there exists $\delta = \delta(\varepsilon )>0$ such that
for any $w,z \in S_X$ with $||w-z|| \leq \delta$, the following estimate holds

$$1 - f_z(w)  = |\ ||z + (w-z)|| - 1 - f_z(w-z) \ | 
\leq \varepsilon ||w-z||/2.$$

So for $x = Rz, y = Rw \in RS_X$ with $||y-x|| \leq {\rm min}\{ R\delta, 2 \}$ it is true that
$$\displaystyle {\rm dist}(y, \Gamma_x)  = R - f_z(y) \leq \varepsilon ||y-x||/2 \leq \varepsilon.$$ 

By assuming $b=2/\delta$ we are done.
 $\black$

\bigskip

The previous fact allows us to say that, in a uniformly smooth space, a "big" sphere intersecting a "small" ball
splits it in two parts that are not "too far" from being slices of the small ball. So Proposition 2.1 leads us to the
following Proposition, which will be crucial for our purposes and may have interest by itself.

\begin{proposition}
\label{fonf}

Let ${\cal B}=\{ B(x_n,R_n) \}_{n=1}^\infty$ be a countable collection of balls in a uniformly smooth
Banach space $X$ such that $R_n$ goes to infinity with $n$.
If $\cal B$ is not locally finite, then it is not point-finite.
\end{proposition}

{\bf Proof.} 
Since ${\cal B}$ is not locally finite, there exist a point $x \in X$ and a sequence $\{ n_k \}$ of integers such that 

\begin{equation}
\label{distance0}
{\rm dist}\bigl( x,B(x_{n_k},R_{n_k} \bigr) \to 0 \ {\rm as} \ k \to \infty.
\end{equation}

Without loss of generality we may assume $x=0$ and $\{ n_k \} = \{ k \}$. 
Moreover, we may assume that $0 \notin B(x_n,R_n)$ for every $n$
and that $ R_n > 2$ and $||x_n|| < R_n + 1$.
For any $n$, let $z_n$ be the point at which the segment $[0,x_n]$ meets $S(x_n,R_n)$ and let
$\Gamma_n$ be the hyperplane supporting $B(x_n,R_n)$ at $z_n$.  Let $0 < \beta_n < 1$ and
$f_n \in S_{X^*}$ be such that $\Gamma_n = \{ t \in X: f_n(t) = \beta_n \}$. Because of (\ref{distance0}),
clearly $||z_n|| \to 0$ and $\beta_n \to 0$ as $ n \to \infty$.

By Fact \ref{F1}, for any $i \in \mathbb N$ big enough
there is $n_i$ such that ${\rm dist}(y,\Gamma_{n_i}) < 1/i$ for every $y \in S(x_{n_i},R_{n_i}) \cap B(z_{n_i},2)$:
hence the set 
$$T_i = \{ t \in B(z_{n_i}, 2): f_{n_i}(t) \geq \beta_{n_i} + 1/i \}$$
is a slice of $B(z_{n_i}, 2)$
contained in $B(x_{n_i},R_{n_i})$. We can choose $n_{i+1} > n_i$ for every $i$. 
Under all our assumptions, for every $i$ big enough (since $B_X \subset B(z_{n_i},2))$ the set
$$S_i = \{ t \in B_X: f_{n_i}(t) \geq \beta_{n_i} + 1/i \}$$ is a slice of $B_X$ that is contained
in $T_i$. 
The sequence $\{ S_i \}$ of such slices of $B_X$ satisfies the assumptions of Proposition
\ref{old}, so it is not point-finite. Since we have 
$$S_i \subset T_i \subset B(x_{n_i},R_{n_i})$$
we are done.  \ $\black$

\bigskip

The following Fact sounds, in some sense, as the converse of
 Fact \ref{F1}. 
 
 \begin{fact}
 \label{F2}
 Let $X$ be uniformly rotund and $b>0$. Let $\{ B_n = B(x_n,R_n) \}_{n=1}^\infty$ be a sequence of balls in $X$
 such that $R_n > b$ and $x_0 \notin {\rm int}B_n$ for every $n$. Put
 $$F_n = \conv \bigl( B_n \setminus {\rm int}B(x_0,b) \bigr)$$
 and let 
 \begin{equation}
 \label{distance}
 {\rm dist}(x_0,F_n) \to 0 \ as \ n \ goes \ to \ \infty.
 \end{equation}
 Then $R_n \to \infty$  with  $n$.
 \end{fact}
 
 {\bf Proof.}
 Let us recall that the following statement is equivalent to $X$ being uniformly rotund: 
 
 \medskip
 
(UR) \ {\it for every $\beta > 0$ 
 there exists $\alpha = \alpha (\beta) > 0$ such that, for any $x \in S_X$ and $f \in S_{X^*}$ with $f(x) = 1$, the slice 
 $\{ t \in B_X: f(t) \geq 1-\alpha \}$ of $B_X$ is contained in $B(x,\beta)$.}
 
 \medskip

Of course we may assume that function $\alpha (\beta) $ is increasing; it goes to 0 as $\beta$ does.

 Without loss of generality 
 may assume $x_0 = 0$.

 Assume to the contrary that some subsequence $\{ R_{n_k} \}$ of $\{ R_n \}$ is bounded: without loss of generality we 
 may assume that $\{ R_n \}$ itself is bounded. For any $n$, put 
$$\varepsilon_n = ||x_n|| - R_n \geq 0 \ \ {\rm and} \ \ 
\tilde{B}_n = B(x_n, R_n + \varepsilon_n);$$ 
moreover, let $f_n \in S_{X *}$ be such that 
 $f_n(\tilde{B}_n) \subset (-\infty,0]$.
 Let $\displaystyle \alpha_n = \alpha ({{b} \over {R_n + \varepsilon_n}})$, where $\alpha = \alpha (\beta)$ is the number defined in condition (UR). We have
 $$\{ t \in \tilde{B}_n: f_n(t) \geq -\alpha_n \} \subset bB_X$$
 
hence
 
 $$\bigl( B_n \setminus {\rm int}bB_X \bigr) \subset \bigl( \tilde{B}_n  \setminus {\rm int}bB_X \bigr) \subset \{ t \in \tilde{B}_n): f_n(t) \leq -\alpha_n \}.$$

 This last set is closed and convex, so
 \begin{equation}
 \label{leq}
 F_n = \conv \bigl( B_n \setminus {\rm int}bB_X \bigr) \subset \{ t \in \tilde{B}_n: f_n(t) \leq -\alpha_n \}.
 \end{equation}
 
 Since $f_n \in S_{X^*}$, (\ref{leq}) implies ${\rm dist}(0,F_n) \geq \alpha_n$ for every $n$: that contradicts (\ref{distance})
because $\{ R_n + \varepsilon_n \}$ is bounded, so $\displaystyle \{ \alpha_n \} = \{ \alpha ({{b} \over {R_n + \varepsilon_n}}) \}$ is far away from 0. 
$\black$

\bigskip

The following three Lemmas are of a technical nature. The first one is something like a local version of
Corson's Theorem \ref{corson}.

\begin{lemma}
\label{lemma1}
Let $X$ be reflexive. Let $x_0 \in X$, $a>b>c>0$ and
 $\cal F$ a  collection of closed convex subsets of $X$ contained in $B(x_0,a)\setminus {\rm int} B(x_0,c)$ such that 
 $\cal F$ covers $B(x_0,a)\setminus {\rm int} B(x_0,b)$. Then $\cal F$ is not locally finite in $X$.
 \end{lemma}
 {\bf Proof.} Without loss of generality we may assume $x_0=0$. Suppose by contradiction 
 that $\cal F$ is locally finite.
 Set ${\cal F}_n =\{ (a/b)^n F : F \in {\cal F}\}$, $n \in \z$ where
 $\alpha F = \{\alpha y : y \in F\}$.  Denote ${\cal F}'=
 \cup_{n \in \z} {\cal F}_n$. Then ${\cal F}'$ covers $X \setminus \{0\}$ and
 $0$ is the only singular point of ${\cal F}'$. Split $X$ into two closed half spaces
 $X=X^+ \cup X^-$ such that $X^+ \cap X^-$ is a hyperplane of $X$, 
 $0 \notin X^+ \cap X^-$ and $0 \in X^-$. Let ${\cal F}^+ =\{ X^+ \cap F : F \in {\cal F}'\}$
 and let ${\cal F}^-$ be the covering of $X^-$ which is a symmetric reflection of
 ${\cal F}^+$ with respect to the hyperplane $X^+\cap X^-$ (made via any line not contained
 in $X^+\cap X^-$). Then ${\cal F}^+ \cup {\cal F}^-$
 is a locally finite covering of $X$ by bounded closed convex sets and we get a contradiction
 with Corson's Theorem \ref{corson}. $\black$
 
\bigskip

 \begin{lemma}
 \label{lemma2}
 Let $X$ be both uniformly rotund and uniformly smooth. Consider a closed hyperplane $X' \subset X$ and  a ball $B(x_0,a),\ x_0 \in X', \ a >0$. Assume that
  ${\cal B}=\{B_n \}_{n=1}^\infty$ is a countable point-finite collection of 
 balls and ${\cal F}=\{F_n \}_{n=1}^\infty$
 is a countable collection of closed convex sets such that 
 $\cal F$ covers $B(x_0,a)\cap X'$,
 $F_n \subset B_n\cap B(x_0,a)$ and
  $x_0 \notin {\rm int} B_n$ for every $n$.
 Then there is a point $y \in B(x_0,a) \cap  X', \ y \neq x_0$, that is a singular point for $\cal F$.
 \end{lemma} 
 {\bf Proof.}
 Assume that every $y \in B(x_0,a)  \cap  X' , \ y \neq x_0$, is a regular point of $\cal F$.
 Take any $b$ such that $a> b>0$ and  define $B'_n =\conv (B_n \setminus {\rm int} B(x_0,b))$,
 $F'_n =B_n'\cap F_n \cap X'$ and ${\cal F}' =\{ F'_n \}_{n=1}^\infty$. 
  We will show
 that there is $c>0$ such that $b>c$ and $F'_n \cap B(x_0,c)=\emptyset$ for every $n$ and then
 our Lemma will follow from Lemma \ref{lemma1} applied to the collection ${\cal F}'$
 with $X$ being replaced by $X'$.
 
 Assume that such $c$ does not exist. Then any arbitrarily
 small neighborhood of $x_0$ intersect $B'_n$ for some $n$ so, since $x_0 \notin {\rm int} B_n$, from Fact \ref{F2} 
we deduce that there is a subsequence of balls in ${\cal B}$ whose radii go to infinity
 for which $x_0$ is a singular point. By Proposition \ref{fonf}, this contradicts 
 the pointwise finiteness of ${\cal B}$. The Lemma is proved. 
 $\black$
 
 \bigskip
  
\begin{lemma}
\label{lemma3}
Let $X$ be both uniformly rotund and uniformly smooth. Let ${\cal B}=\{ B_n \}_{n=1}^\infty = \{ B_n = B(x_n,R_n) \}_{n=1}^\infty$ be a  countable point-finite covering of $X$ by balls.
Put $B^\#_1 =B_1$ and $B^\#_{n+1}= \conv ((B_{n+1} \setminus(B_1 \cup \dots \cup B_n ))$.
Then ${\cal B}^\#= \{ B_n^\# \}_{n=1}^\infty$ is a covering of $X$; moreover, for every $n$ we have that
$B_n^\# \subset B_n$ 
and any $x_0 \in {\rm int} B_n $ is a regular point for ${\cal B}^\#$.

\end{lemma}
{\bf Proof.} Assume that, for some $\tilde n$,  $x_0 \in \Int B_{\tilde n}$
is a singular point for ${\cal B}^\#$. Then for some subsequence 
$\{B_{n_i}^\# \}_{i=1}^\infty$ with  $n_i > \tilde n$ it happens that 
 $ x_0 \notin B_{n_i}$ for every $i$ and $B(x_0, 1/i)$ 
intersects the set $B^\#_{n_j}$  for every $j \geq i$.
Note that $B^\# _{n_i} \subset \conv(B_{n_i} \setminus B_{\tilde n})\subset B_{n_i}$. Then
$B(x_0, 1/i)$ intersects $ \conv(B_{n_j} \setminus B_{\tilde n})$ for $j\geq i$. 
Let $b > 0$ such
that $B(x_0,b) \subset B_{\tilde n}$ (we assumed $x_0 \in {\rm int} B_{\tilde n}$): of course 
$$ \conv(B_{n_j} \setminus B(x_0,b) \supset
 \conv(B_{n_j} \setminus B_{\tilde n})$$ hence, since  and $ x_0 \notin B_{n_i} $ for every $i$, from Fact \ref{F2} we get that 
$ R_{n_i} \lo \infty $ as $i \lo \infty$. By Proposition \ref{fonf} this
contradicts the assumption that $\cal B$ is point-finite. $\black$

 \bigskip
 
 We are now ready to prove our main result.
 
 \bigskip
 
 {\bf Proof of Theorem \ref{goal}.} 
 Assume that ${\cal B} = \{ B_n \}_{n=1}^\infty $ is a point-finite covering by balls of a 
 uniformly rotund and uniformly smooth Banach space $X$.
 Consider the covering ${\cal B}^\#$ from Lemma \ref{lemma3} and let
 $S$ be the set of the points that are singular for ${\cal B}^\#$. By Theorem \ref{corson} we have
 $S \neq \emptyset $ and by Lemma \ref{lemma3} we have  $S\subset \cup_n  \partial B_n$.
 Clearly $S$ is closed in $X$, hence by the Baire cathegory theorem there 
 are $x_0 \in S$,  $a>0$ and $B_m$ such that
 $S \cap B(x_0,a) \subset \partial B_m$. Take a closed hyperplane $X'$ in $X$ passing 
 through $x_0$ and intersecting $B_m$ only at $x_0$. Then, by applying Lemma \ref{lemma2} to the collection 
 ${\cal F}= \{ B^\#_n \cap B(x_0,a) : B^\#_n \in {\cal B}^\#\}$
 of closed convex sets in $X$ with respect to the hyperplane $X'$, we get a contradiction. $\black$

\bigskip
\bigskip
\bigskip

Vladimir P. Fonf 

Department of Mathematics

Ben-Gurion University of the Negev

84105 Beer-Sheva, Israel

E-mail address: fonf@math.bgu.ac.il

\bigskip

Michael Levin 

Department of Mathematics

Ben-Gurion University of the Negev

84105 Beer-Sheva, Israel

E-mail address: mlevine@cs.bgu.ac.il

\bigskip

Clemente Zanco 

Dipartimento di Matematica

Universit\`a degli Studi

Via C. Saldini, 50

20133 Milano MI, Italy

E-mail address: clemente.zanco@unimi.it

ph. ++39 02 503 16164 \ \ \ \ fax  ++39 02 503 16090

\end{document}